\documentclass[12pt]{article}
\usepackage{latexsym,amssymb,amsmath,enumerate,float,geometry,cite}
\newtheorem{theorem}{Theorem}

\newtheorem{lemma}[theorem]{Lemma}

\begin{document}

\title{\Large Bounds on the Exponential Domination Number}
\author{St\'{e}phane Bessy$^1$, Pascal Ochem$^1$, and Dieter Rautenbach$^2$}
\date{}
\maketitle
\vspace{-10mm}
\begin{center}
{\small
$^1$ 
Laboratoire d'Informatique, de Robotique et de Micro\'{e}lectronique de Montpellier (LIRMM),\\
Montpellier, France, \texttt{stephane.bessy@lirmm.fr,pascal.ochem@lirmm.fr}\\[3mm]
$^2$ Institute of Optimization and Operations Research, Ulm University,\\
Ulm, Germany, \texttt{dieter.rautenbach@uni-ulm.de}}
\end{center}

\begin{abstract}
As a natural variant of domination in graphs, 
Dankelmann et al.~[Domination with exponential decay, Discrete Math. 309 (2009) 5877-5883] introduce exponential domination,
where vertices are considered to have some dominating power that decreases exponentially with the distance,
and the dominated vertices have to accumulate a sufficient amount of this power emanating from the dominating vertices.
More precisely, if $S$ is a set of vertices of a graph $G$, 
then $S$ is an exponential dominating set of $G$ if 
$\sum\limits_{v\in S}\left(\frac{1}{2}\right)^{{\rm dist}_{(G,S)}(u,v)-1}\geq 1$
for every vertex $u$ in $V(G)\setminus S$,
where ${\rm dist}_{(G,S)}(u,v)$ is the distance between $u\in V(G)\setminus S$ and $v\in S$ in the graph $G-(S\setminus \{ v\})$.
The exponential domination number $\gamma_e(G)$ of $G$ is the minimum order of an exponential dominating set of $G$.

Dankelmann et al.~show 
$$\frac{1}{4}({\rm d}+2)\leq \gamma_e(G)\leq \frac{2}{5}(n+2)$$
for a connected graph $G$ of order $n$ and diameter ${\rm d}$.
We provide further bounds and in particular strengthen their upper bound.
Specifically, for a connected graph $G$ of order $n$, maximum degree $\Delta$ at least $3$, radius ${\rm r}$ at least $1$, we show
\begin{eqnarray*}
\gamma_e(G) & \geq & \left(\frac{n}{13(\Delta-1)^2}\right)^{\frac{\log_2(\Delta-1)+1}{\log_2^2(\Delta-1)+\log_2(\Delta-1)+1}},\\[3mm]
\gamma_e(G) & \leq & 2^{2{\rm r}-2}\mbox{, and }\\[3mm]
\gamma_e(G) & \leq & \frac{43}{108}(n+2).
\end{eqnarray*}
\end{abstract}

{\small
\begin{tabular}{lp{12.5cm}}
\textbf{Keywords:} & domination, exponential domination\\
\textbf{MSC2010:} & 05C69
\end{tabular}
}

\pagebreak
 
\section{Introduction}

We consider finite, simple, and undirected graphs, and use standard notation and terminology.

A set $D$ of vertices of a graph $G$ is {\it dominating} if every vertex not in $D$ has a neighbor in $D$.
The {\it domination number} $\gamma(G)$ of $G$, defined as the minimum cardinality of a dominating set, 
is one of the most well studied quantities in graph theory \cite{hhs}.
As a natural variant of this classical notion, Dankelmann et al.~\cite{ddems} introduce exponential domination,
where vertices are considered to have some dominating power that decreases exponentially with the distance,
and the dominated vertices have to accumulate a sufficient amount of this power emanating from the dominating vertices.
As a motivation of their model they mention information dissemination within social networks, where
the impact of information decreases every time it is passed on.

Before giving the precise definitions for exponential domination, we mention three closely related well studied notions.
A set $D$ of vertices of a graph $G$ is {\it $k$-dominating} for some positive integer $k$,
if every vertex not in $D$ has at least $k$ neighbors in $D$ \cite{cr,cgs,dghpv,fhv,fj,ha,hv,rv}.
A set $D$ of vertices of a graph $G$ is {\it distance-$k$-dominating} for some positive integer $k$,
if for every vertex not in $D$, there is some vertex in $D$ at distance at most $k$ \cite{adr,bz,chmm,hmv,he,tx}.
Finally, in {\it broadcast domination} \cite{chm,dehhh,e,hl}, each vertex $v$ is assigned an individual dominating power $f(v)$
and dominates all vertices at distance at least $1$ and at most $f(v)$.
Exponential domination shares features with these three notions;
similarly as in $k$-domination, several vertices contribute to the domination of an individual vertex,
similarly as in distance-$k$-domination, vertices dominate others over some distance,
and similarly as in broadcast domination, 
different dominating vertices contribute differently to the domination of an individual vertex depending on the relevant distances.

\medskip

\noindent We proceed to the precise definitions, and also recall some terminology.

Let $G$ be a graph.
The vertex set and the edge set of $G$ are denoted by $V(G)$ and $E(G)$, respectively.
The order $n(G)$ of $G$ is the number of vertices of $G$,
and the size $m(G)$ of $G$ is the number of edges of $G$.
For two vertices $u$ and $v$ of $G$, let ${\rm dist}_G(u,v)$ be the distance in $G$ between $u$ and $v$,
which is the minimum number of edges of a path in $G$ between $u$ and $v$.
If no such path exists, then let ${\rm dist}_G(u,v)=\infty$.
An endvertex is a vertex of degree at most $1$.
For a rooted tree $T$, and a vertex $u$ of $T$, 
let $T_u$ denote the subtree of $T$ rooted in $u$ that contains $u$ as well as all descendants of $u$.
A leaf of a rooted tree is a vertex with no children.
For non-negative integers $d_0,d_1,\ldots,d_k$, 
let $T(d_0,d_1,\ldots,d_k)$ be the rooted tree of depth $k+1$
in which all vertices at distance $i$ from the root have exactly $d_i$ children for every $i$ with $0\leq i\leq k$.
A rooted tree is binary if every vertex has at most two children, 
and a binary tree is full if every vertex other than the leaves has exactly two children.
For a positive integer $k$, let $[k]$ be the set of positive integers at most $k$.

Let $S$ be a set of vertices of $G$.
For two vertices $u$ and $v$ of $G$ with $u\in S$ or $v\in S$,
let ${\rm dist}_{(G,S)}(u,v)$ be the minimum number of edges of a path $P$ in $G$ between $u$ and $v$ 
such that $S$ contains exactly one endvertex of $P$ and no internal vertex of $P$.
If no such path exists, then let ${\rm dist}_{(G,S)}(u,v)=\infty$.
Note that, if $u$ and $v$ are distinct vertices in $S$, then ${\rm dist}_{(G,S)}(u,u)=0$ and ${\rm dist}_{(G,S)}(u,v)=\infty$.

For a vertex $u$ of $G$, let 
$$w_{(G,S)}(u)=\sum\limits_{v\in S}\left(\frac{1}{2}\right)^{{\rm dist}_{(G,S)}(u,v)-1},$$
where $\left(\frac{1}{2}\right)^{\infty}=0$.
Note that $w_{(G,S)}(u)=2$ for $u\in S$.

If $w_{(G,S)}(u)\geq 1$ for every vertex $u$ of $G$, then $S$ is an {\it exponential dominating set} of $G$.
The {\it exponential domination number $\gamma_e(G)$} is the minimum order of an exponential dominating set of $G$,
and an exponential dominating set of $G$ of order $\gamma_e(G)$ is {\it minimum}.
By definition, every dominating set is also an exponential dominating set, 
which implies $\gamma_e(G)\leq \gamma(G)$ for every graph $G$.

\medskip

\noindent The following summarizes the main results of Dankelmann et al.~\cite{ddems}.

\begin{theorem}[Dankelmann et al.~\cite{ddems}]\label{theoremd}
If $G$ is a connected graph of diameter ${\rm diam}(G)$, then 
$$\frac{1}{4}({\rm diam}(G)+2)\leq \gamma_e(G)\leq \frac{2}{5}(n(G)+2).$$
\end{theorem}
Dankelmann et al.~\cite{ddems} discuss the tightness of their bounds.
They show that the lower bound is satisfied with equality for the path $P_n$ of order $n$ with $n\equiv 2\mod 4$,
and they construct a sequence of trees $T$ for which $\frac{\gamma_e(T)}{n(T)+2}$ tends to $\frac{3}{8}$.
Finally, they describe one specific tree $T$ with $\frac{\gamma_e(T)}{n(T)+2}=\frac{144}{377}\approx 0.382$,
and ask whether there are trees $T$ with $\frac{144}{377}<\frac{\gamma_e(T)}{n(T)+2}\leq \frac{2}{5}$.

Note that the lower bound in Theorem \ref{theoremd} implies $\gamma_e(G)=\Omega(\log n(G))$ for graphs $G$ of bounded maximum degree,
because the diameter of such graphs is $\Omega(\log n(G))$.
Our first result is a polynomial, and not just logarithmic, lower bound.

\begin{theorem}\label{theoremlbgen}
If $G$ is a graph of maximum degree $\Delta(G)$ at least $3$, then 
$$\gamma_e(G)\geq \left(\frac{n(G)}{13(\Delta(G)-1)^2}\right)^{\frac{\log_2(\Delta(G)-1)+1}{\log_2^2(\Delta(G)-1)+\log_2(\Delta(G)-1)+1}}.$$
\end{theorem}
As our second result, we show that $\gamma_e(G)$ is not only lower bounded but in fact also upper bounded in terms of the diameter of $G$,
or rather the radius of $G$.

\begin{theorem}\label{theoremdiam}
If $G$ is a connected graph of radius ${\rm rad}(G)$ at least $1$, then $\gamma_e(G)\leq 2^{2{\rm rad}(G)-2}$.
\end{theorem}
Surprisingly, the bound in Theorem \ref{theoremdiam} is tight as we show by constructing a suitable example.

As our third result, we improve the upper bound in Theorem \ref{theoremd} as follows.

\begin{theorem}\label{theorem1}
If $G$ is a connected graph, then $\gamma_e(G)\leq \frac{43}{108}(n(G)+2)$.
\end{theorem}
Note that $\frac{43}{108}\approx 0.398$.

All proofs and further discussion are postponed to the next section.

\section{Proofs}

\noindent {\it Proof of Theorem \ref{theoremlbgen}:} 
Let $\Delta=\Delta(G)$, and $\alpha=1-\frac{1}{\log_2(\Delta-1)+1}$.
Let $S$ be an exponential dominating set of $G$.
Let $H$ arise from $G$ by removing all edges between vertices in $S$.
Clearly, $S$ is still an exponential dominating set of $H$.

Let $k=|S|$.

Let 
\begin{eqnarray*}
A & = & \left\{ v\in V(G)\setminus S: {\rm dist}_H(v,S)\leq \alpha\log_2(k)\right\}\mbox{ and}\\
B & = & \left\{ v\in V(G)\setminus S: {\rm dist}_H(v,S)> \alpha\log_2(k)\right\},
\end{eqnarray*}
where ${\rm dist}_H(v,S)=\min\{ {\rm dist}_H(v,u):u\in S\}$.

For $u\in S$, let
\begin{eqnarray*}
C(u) & = & \left\{ v\in B: {\rm dist}_H(v,u)\leq \log_2(k)+2\right\}.
\end{eqnarray*}
Since in a graph of maximum degree $\Delta$, 
there are at most $\frac{\Delta}{\Delta-2}\left((\Delta-1)^d-1\right)\leq 3\left((\Delta-1)^d-1\right)$
vertices at distance between $1$ and $d$ from any given vertex, we obtain 
$$|A|\leq 3k\left((\Delta-1)^{\alpha\log_2(k)}-1\right)=
3k^{\frac{\log_2^2(\Delta-1)+\log_2(\Delta-1)+1}{\log_2(\Delta-1)+1}}-3k.$$
Let 
$${\cal R}=\left\{ (u,v):u\in S,v\in C(u)\right\}.$$
Since $(u,v)\in {\cal R}$ implies ${\rm dist}_H(v,u)\leq \log_2(k)+2$,
we obtain that, for every $u$ in $S$, 
there are at most 
$3\left((\Delta-1)^{\log_2(k)+2}-1\right)\leq 3(\Delta-1)^2k^{\log_2(\Delta-1)}$ vertices $v$ with $(u,v)\in {\cal R}$,
which implies 
$$|{\cal R}|\leq 3k(\Delta-1)^2k^{\log_2(\Delta-1)}.$$
If there is some $v$ in $B$ such that there are less than $\frac{1}{4}k^{\alpha}$ vertices $u$ with $(u,v)\in {\cal R}$,
then ${\rm dist}_H(v,u')>\log_2(k)+2$ for more than $k-\frac{1}{4}k^{\alpha}$ vertices $u'$ in $S$.
Since $v\in B$ implies ${\rm dist}_H(v,S)>\alpha\log_2(k)$, we obtain
$$w_{(H,S)}(v)<\frac{1}{4}k^{\alpha}\left(\frac{1}{2}\right)^{\alpha\log_2(k)-1}
+\left(k-\frac{1}{4}k^{\alpha}\right)\left(\frac{1}{2}\right)^{\left(\log_2(k)+2\right)-1}
\leq \frac{1}{2}+\frac{1}{2}\frac{\left(k-\frac{1}{4}k^{\alpha}\right)}{k}<1,$$
which is a contradiction.
Hence, for every $v$ in $B$, there are at least $\frac{1}{4}k^{\alpha}$ vertices $u$ with $(u,v)\in {\cal R}$,
which implies 
$$|{\cal R}|\geq \frac{1}{4}k^{\alpha}|B|.$$
Combining the upper and the lower bound on $|{\cal R}|$, we obtain 
$$|B|\leq 12(\Delta-1)^2k^{\log_2(\Delta-1)+1-\alpha}=12(\Delta-1)^2 k^{\frac{\log_2^2(\Delta-1)+\log_2(\Delta-1)+1}{\log_2(\Delta-1)+1}}.$$
Altogether, we obtain
\begin{eqnarray*}
n(G)&= & |S|+|A|+|B|\\
&\leq & k+3k^{\frac{\log_2^2(\Delta-1)+\log_2(\Delta-1)+1}{\log_2(\Delta-1)+1}}-3k
+12(\Delta-1)^2 k^{\frac{\log_2^2(\Delta-1)+\log_2(\Delta-1)+1}{\log_2(\Delta-1)+1}}\\
&\leq & 13(\Delta-1)^2 k^{\frac{\log_2^2(\Delta-1)+\log_2(\Delta-1)+1}{\log_2(\Delta-1)+1}},
\end{eqnarray*}
which implies the desired bound.
$\Box$

\medskip

\noindent It is not difficult to improve the constant $13$ in Theorem \ref{theoremlbgen} by adding some technicalities.
For the sake of simplicity, we decided not to do so.

\medskip

\noindent {\it Proof of Theorem \ref{theoremdiam}:}
Since $G$ has a rooted spanning tree $T$ of depth at most ${\rm rad}(G)$, and $\gamma_e(G)\leq \gamma_e(T)$, 
it suffices to show $\gamma_e(T)\leq 2^{2d-2}$ for a rooted tree $T$ of depth $d$ at least $1$.
The proof is by induction on the depth $d$ of $T$.

If $d=1$, then the root $r$ of $T$ forms an exponential dominating set of $T$, and hence, $\gamma_e(T)=1=2^{2\cdot 1-2}$.
If $d=2$, then four children of $r$ form an exponential dominating set of $T$,
and, if $r$ does not have four children, then the set of all its children forms an exponential dominating set of $T$.
Hence, $\gamma_e(T)\leq 4=2^{2\cdot 2-2}$,
and we may assume that $d\geq 3$.

If $S$ is a set of $2^{2d-2}$ parents of leaves of $T$, 
then, since the distance between any vertex in $S$ and any other vertex of $T$ is at most $2d-1$,
and $\left(\frac{1}{2}\right)^{(2d-1)-1}|S|=1$,
the set $S$ is an exponential dominating set of $T$.
Hence, we may assume that the set $S_0$ of all parents of leaves of $T$ has less than $2^{2d-2}$ elements,
and is not an exponential dominating set of $T$.

Suppose that $S_0$ has at least $\frac{1}{8}\cdot 2^{2d-2}$ elements.
Let $u$ be any vertex of $T$.
If $u$ has depth at least $d-2$, then some vertex in $S_0$ has distance at most $1$ from $u$,
which implies $w_{(T,S_0)}(u)\geq 1$.
If $u$ has depth at most $d-3$, then the distance between $u$ and any vertex in $S_0$ is at most $2d-4$.
Since $\left(\frac{1}{2}\right)^{(2d-4)-1}|S_0|\geq 1$, we obtain $w_{(T,S_0)}(u)\geq 1$ also in this case,
which implies the contradiction that $S_0$ is an exponential dominating set of $T$.
Hence, $S_0$ has less than $\frac{1}{8}\cdot 2^{2d-2}$ elements.

For a vertex $u$ of $T$, let $w(u)=w_{(T_u,S_0\cap V(T_u))}(u)$.
Let $T_1$ arise from $T$ by removing every vertex $u$ such that $w(v)\geq 1$ for every vertex $v$ in $V(T_u)$.
By the choice of $S_0$, this construction implies that $T_1$ is a rooted tree of depth at most $d-3$.
Note that $T_1$ might have depth $0$, that is, it may consist only of the root.
By induction, $T_1$ has an exponential dominating set $S_1$ of order at most $\max\left\{1,2^{2(d-3)-2}\right\}$.
Now, $S_0\cup S_1$ is an exponential dominating set of $T$ of order at most 
$\frac{1}{8}\cdot 2^{2d-2}+\max\left\{1,2^{2(d-3)-2}\right\}\leq 2^{2d-2}$,
which completes the proof.
$\Box$

\medskip

\noindent In order to show that the bound in Theorem \ref{theoremdiam} is tight,
we need a simple observation concerning binary trees.

\begin{lemma}\label{lemmabinarytree}
If $T$ is a binary tree with root $r$, and $S$ is a set of vertices of $T$, then $w_{(T,S)}(r)\leq 2$
with equality if and only if $T$ contains a full binary subtree $F$ with root $r$ such that $V(F)\cap S$ is the set of leaves of $F$.
\end{lemma}
{\it Proof:} The proof is by induction on the depth $d$ of $T$.
If $d=0$ or $r\in S$, then the statement is trivial.
Hence, we may assume that $d\geq 1$ and $r\not\in S$.
Let $r_1,\ldots,r_k$ for some $k\in [2]$ be the children of $r$.
For $i\in [k]$, let $T_i$ be the subtree of $T$ rooted in $r_i$, and let $S_i=S\cap V(T_i)$.
Since $w_{(T,S)}(r)=\frac{1}{2}\sum_{i=1}^k w_{(T_i,S_i)}(r_i)$, we obtain, by induction, 
that $w_{(T,S)}(r)\leq 2$ with equality if and only if $k=2$, and $w_{(T_1,S_1)}(r_1)=w_{(T_2,S_2)}(r_2)=2$.
Now, $w_{(T_1,S_1)}(r_1)=w_{(T_2,S_2)}(r_2)=2$ is equivalent with the existence of suitable full binary subtrees $F_1$ and $F_2$ as described in the statement. 
Since the existence of $F_1$ and $F_2$ is clearly equivalent with the existence of the subtree $F$ as described in the statement,
the proof is complete. $\Box$

\medskip

\noindent For some positive integer $d$, 
let the rooted tree $T$ arise by attaching $2^{2d-2}$ disjoint copies of the full binary tree 
$$T(\underbrace{2,\ldots,2,0}_{d})$$ 
of depth $d-1$ to the root $r$ of $T$.
By Theorem \ref{theoremdiam}, we have $\gamma_e(T)\leq 2^{2d-2}$.
In fact, we are going to show that $\gamma_e(T)=2^{2d-2}$.
Therefore, let $S$ be a minimum exponential dominating set of $T$ that does not contain any leaf (notice that
if $S$ contains a leaf, then we could replace this leaf by its parent and still have an exponential dominating set).
Assume that $S$ contains a vertex $u$ that is neither the root $r$ nor a parent of a leaf.
If $S$ does not contain the parent $v$ of some leaf $v'$ of $T_u$, then, as $v'$ must be dominated, we must have $w_{(T,S)}(v)\ge 2$.
By Lemma~\ref{lemmabinarytree}, this implies in particular that the second child of $v$, which is a leaf,
is in $S$, contradicting the fact that $S$ does not contain any leaf of $T$. So $S$ contains all the parents
of the leaves of $T_u$, and we have $w_{(T_u,S\setminus \{u\})}(u)=2$. 
Now, $S\setminus \{ u\}$ is an exponential dominating set
of $T$, a contradiction.
Hence, $S\setminus \{ r\}$ contains only parents of leaves. 
Suppose that $|S|<2^{2d-2}$. 
This implies that $d\geq 2$, and that there is some child $x$ of the root $r$ such that $S\cap V(T_x)=\emptyset$.
Since $S$ is an exponential dominating set, it follows in particular that $S$ does not contain $r$.
So we have ${\rm dist}_T(u,v)=2d-1$ for every vertex $u$ in $S$ and every leaf $v$ in $V(T_x)$, and
we obtain $\left(\frac{1}{2}\right)^{(2d-1)-1}|S|\geq 1$,
that is, $|S|\geq 2^{2d-2}$, which is a contradiction.
Altogether, it follows that $\gamma_e(T)=2^{2d-2}$.

\medskip

\noindent Our next goal is to prove Theorem \ref{theorem1}.

Similarly as the proof of Theorem \ref{theoremd} in \cite{ddems}, the
proof of Theorem \ref{theorem1} is based on an inductive argument that
uses local reductions.  Unfortunately, the non-local character of
exponential domination makes it unlikely that a local approach can
lead to a best-possible result.  Even in order to achieve a very small
improvement of the upper bound in Theorem \ref{theoremd}, the approach
makes it necessary to consider a large number of cases and specific
configurations.  Since our goal was rather to obtain a constant lower
than $2/5$ for the upper bound in Theorem \ref{theoremd} than to obtain
the best-possible result, we tried to limit the number of cases as
much as possible for the sake of simplicity.  There are several parts
of our proof though, where further obvious improvements are possible
at the cost of considering more cases.

In the next subsection we collect several auxiliary results,
and in Subsection \ref{sec2.2} we prove Theorem \ref{theorem1}.

\subsection{Auxiliary results}\label{sec2.1}

\begin{lemma}\label{lemma1}
Let $T$ be a tree.
\begin{enumerate}[(i)]
\item If ${\rm diam}(T)\leq 2$, then $\gamma_e(T)=1$.
\item If ${\rm diam}(T)=3$, then $\gamma_e(T)=2$ and $n(T)\geq 4$.
\item If ${\rm diam}(T)=4$, then let $u$ be the central vertex of $T$.
Let $u$ have $k$ neighbors that are endvertices and $\ell\geq 2$ neighbors that are not endvertices.
\begin{enumerate}[(a)]
\item If $\ell=2$, then $\gamma_e(T)=2$ and $n(T)\geq 5$.
\item If $\ell=3$, then $\gamma_e(T)\leq 3$ and $n(T)\geq 7$.
\item If $\ell\geq 4$, then $\gamma_e(T)\leq 4$ and $n(T)\geq 9$.
\end{enumerate}
\end{enumerate}
\end{lemma}
{\it Proof:} Since the proofs of (i) and (ii) are straightforward, we
only give details for the proof of (iii).  Since $T$ has diameter $4$,
no single vertex forms an exponential dominating set of $T$, which
implies $\gamma_e(T)>1$.  If $\ell\leq 3$, then the neighbors of $u$
that are no endvertices form an exponential dominating set of $T$,
which implies $\gamma_e(T)\leq \ell$.  If $\ell\geq 4$, then four
neighbors of $u$ that are no endvertices form an exponential
dominating set of $T$, which implies $\gamma_e(T)\leq 4$.  Since
$n(T)\geq 1+k+2\ell$, the lower bounds on the order of $T$ follow.
$\Box$

\medskip

\noindent For the rest of this subsection, let $T$ be a tree of diameter at least $5$.
We root $T$ in a vertex of maximum eccentricity, that is, the depth of $T$ is at least $3$.

\begin{lemma}\label{lemma2}
Let $u$ be a vertex of $T$, and let $v_1,\ldots,v_k$ be some children of $u$.\\
If one of the following conditions (\ref{r1}) to (\ref{r16}) holds,
then there is a tree $T'$ with $n(T')<n(T)$ and $\gamma_e(T)\leq \gamma_e(T')+\frac{5}{13}(n(T)-n(T'))$.
\begin{enumerate}[(i)]
\item\label{r1} $k=2$, and $v_1$ and $v_2$ are leaves.
\item\label{r2} $k=1$, and $T_{v_1}\cong T(1,1,0)$.
\item\label{r3} $k=2$, $v_1$ is a leaf, and $T_{v_2}\cong T(1,0)$.
\item\label{r4} $k=4$, and $T_{v_i}\cong T(1,0)$ for $i\in [4]$.
\item\label{r5} $k=2$, $v_1$ is a leaf, and $T_{v_2}\cong T(2,1,0)$.
\item\label{r6} $k=2$, $v_1$ is a leaf, and $T_{v_2}\cong T(3,1,0)$.
\item\label{r7} $k=6$, and $T_{v_i}\cong T(2,1,0)$ for $i\in [6]$.
\item\label{r8} $k=3$, and $T_{v_i}\cong T(3,1,0)$ for $i\in [3]$.
\item\label{r12} $k=4$, $T_{v_1}\cong T(3,1,0)$, and $T_{v_i}\cong T(2,1,0)$ for $i\in [4]\setminus [1]$.
\item\label{r13} $k=4$, $T_{v_i}\cong T(3,1,0)$ for $i\in [2]$, and $T_{v_i}\cong T(2,1,0)$ for $i\in [4]\setminus [2]$.
\item\label{r9} $k=4$, $T_{v_1}\cong T(1,0)$, and $T_{v_i}\cong T(2,1,0)$ for $i\in [4]\setminus [1]$.
\item\label{r10} $k=3$, $T_{v_1}\cong T(1,0)$, and $T_{v_i}\cong T(3,1,0)$ for $i\in [3]\setminus [1]$.
\item\label{r11} $k=4$, $T_{v_1}\cong T(1,0)$, $T_{v_2}\cong T(3,1,0)$, and $T_{v_i}\cong T(2,1,0)$ for $i\in [4]\setminus [2]$.
\item\label{r14} $k=4$, $T_{v_i}\cong T(1,0)$ for $i\in [3]$, and $T_{v_4}\cong T(3,1,0)$.
\item\label{r17} $k=4$, $T_{v_i}\cong T(1,0)$ for $i\in [3]$, and $T_{v_4}\cong T(2,1,0)$.
\item\label{r15} $k=4$, $T_{v_i}\cong T(1,0)$ for $i\in [2]$, and $T_{v_i}\cong T(2,1,0)$ for $i\in [4]\setminus [2]$.
\item\label{r16} $k=4$, $T_{v_i}\cong T(1,0)$ for $i\in [2]$, $T_{v_3}\cong T(3,1,0)$, and $T_{v_4}\cong T(2,1,0)$.
\end{enumerate}
\end{lemma}
{\it Proof:} We consider different cases corresponding to the above conditions.
In each case, we construct a suitable tree $T'$ with $n(T')<n(T)$.
Throughout the proof, let $S'$ be a minimum exponential dominating set of $T'$.

If (\ref{r1}) occurs, then let $T'=T-v_2$.
If $u\in S'$ or $u,v_1\not\in S'$, then $S'$ is also an exponential dominating set of $T$.
If $u\not\in S'$ and $v_1\in S'$, then $S=(S'\setminus \{ v_1\})\cup \{ u\}$ is an exponential dominating set of $T$.
We obtain $\gamma_e(T)\leq \gamma_e(T')+0(n(T)-n(T'))$.

If (\ref{r2}) occurs, then let $T'=T-V(T_{v_1})$.
If $w$ is the neighbor of $v_1$ in $T_{v_1}$, then $S'\cup \{ w\}$ is an exponential dominating set of $T$.
We obtain $\gamma_e(T)\leq \gamma_e(T')+\frac{1}{3}(n(T)-n(T'))$.

\begin{figure}[H]
\begin{center}
\unitlength 1mm 
\linethickness{0.4pt}
\ifx\plotpoint\undefined\newsavebox{\plotpoint}\fi 
\begin{picture}(55,33)(0,0)
\put(0,10){\circle*{2}}
\put(10,10){\circle*{2}}
\put(20,10){\circle*{2}}
\put(0,20){\circle*{2}}
\put(10,20){\circle*{2}}
\put(20,20){\circle*{2}}
\put(40,20){\circle*{2}}
\put(40,10){\circle*{2}}
\put(40,0){\circle*{2}}
\put(30,10){\circle*{2}}
\put(30,0){\circle*{2}}
\put(50,10){\circle*{2}}
\put(50,0){\circle*{2}}
\put(20,30){\circle*{2}}
\put(50,0){\line(0,1){10}}
\put(50,10){\line(-1,1){10}}
\put(40,20){\line(0,-1){20}}
\put(30,0){\line(0,1){10}}
\put(30,10){\line(1,1){10}}
\put(40,20){\line(-2,1){20}}
\put(20,30){\line(0,-1){20}}
\put(20,30){\line(-1,-1){10}}
\put(10,20){\line(0,-1){10}}
\put(20,30){\line(-2,-1){20}}
\put(0,20){\line(0,-1){10}}
\put(8,22){\line(0,-1){4}}
\put(8,18){\line(1,0){1}}
\put(9,18){\line(1,0){16}}
\put(25,18){\line(0,-1){10}}
\put(25,8){\line(1,0){27}}
\put(52,8){\line(0,1){2}}
\put(52,10){\line(0,1){2}}
\put(52,12){\line(-1,0){24}}
\put(28,12){\line(0,1){0}}
\put(28,12){\line(0,1){10}}
\put(28,22){\line(-1,0){20}}
\put(55,10){\makebox(0,0)[cc]{$W$}}
\put(20,33){\makebox(0,0)[cc]{$u$}}
\put(40,23){\makebox(0,0)[cc]{$v_4$}}
\put(0,23){\makebox(0,0)[cc]{$w$}}
\end{picture}
\end{center}
\caption{The configuration in case (\ref{r14}).}\label{figr14}
\end{figure}
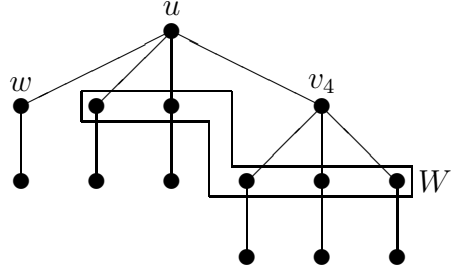

For the remaining cases, let $T'=T-\bigcup_{i=1}^kV(T_{v_i})$.
We specify a vertex $w$ and a set $W$ with the following properties:
\begin{itemize}
\item If $u\not\in S'$, then $S'\cup W$ is an exponential dominating set of $T$.
\item If $u\in S'$, then $(S'\setminus \{ u\})\cup \{ w\}\cup W$ is an exponential dominating set of $T$.
\item $\gamma_e(T)\leq \gamma_e(T')+c(n-n')$ with $c \le \frac{5}{13}$.
\end{itemize}
We leave it to the reader to verify the straightforward details. 
\begin{itemize}
\item (\ref{r3}):
Let $w=v_1$ and $W=\{ v_2\}$. We obtain $c=\frac{1}{3}$.
\item (\ref{r4}):
Let $w=v_1$ and $W=\{ v_2,v_3,v_4\}$.  We obtain $c=\frac{3}{8}$.
\item (\ref{r5}) or (\ref{r6}):
Let $w=v_1$ and let $W$ be the set of children of $v_2$.  We respectively obtain $c=\frac{1}{3}$ and $c=\frac{3}{8}$.
\item (\ref{r7}) or (\ref{r8}) or (\ref{r12}) or (\ref{r13}):
Let $w$ be a child of $v_1$ and let $W$ be the set of children of $v_1,\ldots,v_k$ except for $w$.
We respectively obtain $c=\frac{11}{30}$,  $c=\frac{8}{21}$, $c=\frac{8}{22}$, and $c=\frac{9}{24}$.
\item (\ref{r9}) or (\ref{r10}) or (\ref{r11}):
Let $w$ be $v_1$ and let $W$ be the set of children of $v_2,\ldots,v_k$.
We respectively obtain  $c=\frac{6}{17}$, $c=\frac{6}{16}$, and $c=\frac{7}{19}$.
\item (\ref{r14}) or (\ref{r17}):
Let $w=v_1$ and let $W$ be the set containing $v_2$, $v_3$ as well as the children of $v_4$.
We respectively obtain $c=\frac{5}{13}$ and $c=\frac{4}{11}$.
See Figure \ref{figr14} for an illustration of case (\ref{r14}).
\item (\ref{r15}) or (\ref{r16}): 
Let $w=v_1$ and let $W$ be the set containing $v_2$ as well as the children of $v_3$ and $v_4$.
We respectively obtain $c=\frac{5}{14}$ and $c=\frac{6}{16}$.
\end{itemize}
Note that the factor $\frac{5}{13}$ comes from case (\ref{r14}).
The other cases actually lead to smaller factors.
$\Box$

\begin{lemma}\label{lemma3}
Let $w$ be a vertex of $T$, and let $X$ be a set of some children of
$w$ such that $T_x$ has depth at most $2$ for every $x$ in $X$, and
$T_x$ has depth $2$ for at least one $x$ in $X$.\\ If there is no tree
$T'$ with $n(T')<n(T)$ and $\gamma_e(T)\leq
\gamma_e(T')+\frac{7}{18}(n(T)-n(T'))$, then there are non-negative
integers $k_1$, $k_2$, and $k_3$ such that
\begin{itemize}
\item $k_1+k_2+k_3=|X|$.
\item $T_x\cong T(1,0)$ for $k_1$ vertices $x$ in $X$. 
\item $T_x\cong T(2,1,0)$ for $k_2$ vertices $x$ in $X$.
\item $T_x\cong T(3,1,0)$ for $k_3$ vertices $x$ in $X$.
\item Furthermore, $k_1$, $k_2$, and $k_3$ satisfy the following restrictions.
\begin{enumerate}[(a)]
\item\label{lemma8a} $k_3\leq 2$, and if $k_3=2$, then $(k_1,k_2,k_3)=(0,0,2)$.
\item\label{lemma8b} If $k_3=1$, then $k_2\leq 2$.
\item\label{lemma8c} If $k_1\geq 1$ and $k_3=1$, then $k_1\leq 2$ and $k_2\leq 1$.
\item\label{lemma8d} If $k_2=1$ and $k_3=1$, then $k_1\leq 1$.
\item\label{lemma8e} If $k_3=0$, then $k_2\leq 5$.
\item\label{lemma8f} If $k_1\geq 1$ and $k_3=0$, then $k_2\leq 2$.
\item\label{lemma8g} If $k_1\geq 1$, $k_2=2$, and $k_3=0$, then $k_1=1$.
\item\label{lemma8h} If $k_1\geq 1$, $k_2=1$, and $k_3=0$, then $k_1\leq 2$.
\end{enumerate}
\end{itemize}
\end{lemma}
{\it Proof:} Since $\frac{5}{13}<\frac{7}{18}$,
we may assume, by Lemma \ref{lemma2}, that $T$ does not contain any of the substructures described in that lemma.
By Lemma \ref{lemma2}(\ref{r1}), $T_x\cong T(1,0)$ for every $x$ in $X$ such that $T_x$ has depth $1$.
Let $k_1$ be the number of $x$ in $X$ such that $T_x\cong T(1,0)$.
By Lemma \ref{lemma2}(\ref{r1}) to (\ref{r4}), 
$T_x\cong T(2,1,0)$ or $T_x\cong T(3,1,0)$ for every $x$ in $X$ such that $T_x$ has depth $2$.
Let $k_2$ and $k_3$ be the numbers of $x$ in $X$ such that $T_x\cong T(2,1,0)$ and $T_x\cong T(3,1,0)$, respectively.
Since $T_x$ has depth $2$ for at least one $x$ in $X$, we have $k_2+k_3\geq 1$.
By Lemma \ref{lemma2}(\ref{r3}), (\ref{r5}), and (\ref{r6}), $T_x$ has depth at least $1$ for every $x$ in $X$,
which implies $k_1+k_2+k_3=|X|$.
By Lemma \ref{lemma2}(\ref{r8}), we have $k_3\leq 2$.

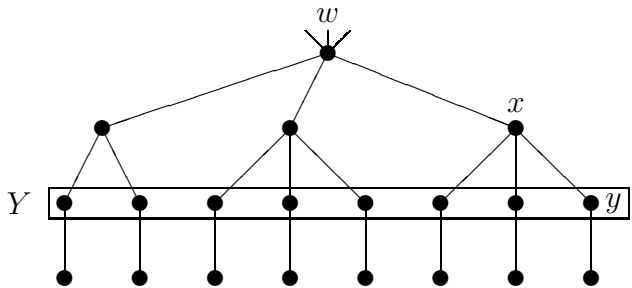
\begin{figure}[H]
\begin{center}
\unitlength 1mm 
\linethickness{0.4pt}
\ifx\plotpoint\undefined\newsavebox{\plotpoint}\fi 
\begin{picture}(85,35)(0,0)
\put(10,0){\circle*{2}}
\put(20,0){\circle*{2}}
\put(30,0){\circle*{2}}
\put(40,0){\circle*{2}}
\put(50,0){\circle*{2}}
\put(60,0){\circle*{2}}
\put(70,0){\circle*{2}}
\put(80,0){\circle*{2}}
\put(10,10){\circle*{2}}
\put(20,10){\circle*{2}}
\put(30,10){\circle*{2}}
\put(40,10){\circle*{2}}
\put(50,10){\circle*{2}}
\put(60,10){\circle*{2}}
\put(70,10){\circle*{2}}
\put(80,10){\circle*{2}}
\put(10,10){\line(0,-1){10}}
\put(20,10){\line(0,-1){10}}
\put(30,10){\line(0,-1){10}}
\put(40,10){\line(0,-1){10}}
\put(50,10){\line(0,-1){10}}
\put(60,10){\line(0,-1){10}}
\put(70,10){\line(0,-1){10}}
\put(80,10){\line(0,-1){10}}
\put(15,20){\circle*{2}}
\put(40,20){\circle*{2}}
\put(70,20){\circle*{2}}
\put(80,10){\line(-1,1){10}}
\put(70,20){\line(0,-1){10}}
\put(50,10){\line(-1,1){10}}
\put(40,20){\line(0,-1){10}}
\put(20,10){\line(-1,2){5}}
\put(15,20){\line(-1,-2){5}}
\put(45,30){\circle*{2}}
\put(30,10){\line(1,1){10}}
\put(40,20){\line(1,2){5}}
\put(45,30){\line(5,-2){25}}
\put(70,20){\line(-1,-1){10}}
\put(45,30){\line(-3,-1){30}}
\multiput(45,30)(.03370787,.03370787){89}{\line(0,1){.03370787}}
\put(45,30){\line(0,1){3}}
\multiput(45,30)(-.03370787,.03370787){89}{\line(0,1){.03370787}}
\put(42,33){\line(0,1){0}}
\put(45,35){\makebox(0,0)[cc]{$w$}}
\put(70,23){\makebox(0,0)[cc]{$x$}}
\put(83,10){\makebox(0,0)[cc]{$y$}}
\put(8,8){\framebox(77,4)[cc]{}}
\put(4,10){\makebox(0,0)[cc]{$Y$}}
\end{picture}
\end{center}
\caption{$(k_1,k_2,k_3)=(0,1,2)$.}\label{figlemma8}
\end{figure}

Suppose now that $k_3=2$.
By Lemma \ref{lemma2}(\ref{r13}) and (\ref{r10}), $k_1=0$ and $k_2\leq 1$,
which implies $(k_1,k_2,k_3)\in \{ (0,0,2),(0,1,2)\}$.
If $(k_1,k_2,k_3)=(0,1,2)$, then let $T'$ arise from $T$ by removing all descendants of $w$ except for one child $x$ of $w$.
Let $S'$ be a minimum exponential dominating set of $T'$.
Clearly, we may assume that $x\not\in S'$.
Let $Y$ be the set of the eight descendants of $w$ at distance $2$ from $w$.
Let $y$ be a vertex in $Y$ with a neighbor of degree $4$.
See Figure \ref{figlemma8} for an illustration.

\noindent If $w\not\in S'$, then $w_{(T',S')}(x)\geq 1$ implies $w_{(T',S')}(w)\geq 2$, and hence,
$S'\cup (Y\setminus \{ y\})$ is an exponential dominating set of $T$.
If $w\in S'$, then $(S'\setminus \{ w\})\cup Y$ is an exponential dominating set of $T$.
In both cases, $\gamma_e(T)\leq \gamma_e(T')+\frac{7}{18}(n(T)-n(T'))$.
Hence, $k_3=2$ implies $(k_1,k_2,k_3)=(0,0,2)$.

If $k_3=1$, then Lemma \ref{lemma2}(\ref{r12}) implies $k_2\leq 2$.
If $k_1\geq 1$ and $k_3=1$, then Lemma \ref{lemma2}(\ref{r11}) and (\ref{r14}) imply $k_2\leq 1$ and $k_1\leq 2$.
If $k_2=1$ and $k_3=1$, then Lemma \ref{lemma2}(\ref{r16}) implies $k_1\leq 1$.
If $k_3=0$, then Lemma \ref{lemma2}(\ref{r7}) implies $k_2\leq 5$.
If $k_1\geq 1$ and $k_3=0$, then Lemma \ref{lemma2}(\ref{r9}) implies $k_2\leq 2$.
If $k_1\geq 1$, $k_2=2$, and $k_3=0$, then Lemma \ref{lemma2}(\ref{r15}) implies $k_1=1$.
If $k_1\geq 1$, $k_2=1$, and $k_3=0$, then Lemma \ref{lemma2}(\ref{r17}) implies $k_1\leq 2$.
$\Box$

\medskip

\noindent Note that in Lemma \ref{lemma2} and Lemma \ref{lemma3} we
consider only some and not necessarily all children of $u$ and $w$,
respectively.

A vertex $w$ of $T$ has {\it type $(k_1,k_2,k_3)$} for non-negative integers $k_1$, $k_2$, and $k_3$ with $k_2+k_3\geq 1$,
if 
\begin{itemize}
\item $k_1$, $k_2$, and $k_3$ satisfy the restrictions stated in Lemma \ref{lemma3}(\ref{lemma8a}) to (\ref{lemma8h}),
\item $w$ has exactly $k_1+k_2+k_3$ children,
\item $T_x\cong T(1,0)$ for $k_1$ children $x$ of $w$, 
\item $T_x\cong T(2,1,0)$ for $k_2$ children $x$ of $w$, and
\item $T_x\cong T(3,1,0)$ for $k_3$ children $x$ of $w$.
\end{itemize}
Note that if $w$ has some type, then  $T_w$ has depth $3$.

\begin{lemma}\label{lemma4}
Let the vertex $w$ of $T$ have type $(k_1,k_2,k_3)$.\\
If $(k_1,k_2,k_3)\not\in \{ (0,0,2),(1,0,1),(2,0,1),(2,1,0)\}$, 
then there is a tree $T'$ with $n(T')\leq n(T)-6$ and $\gamma_e(T)\leq \gamma_e(T')+\frac{7}{18}(n(T)-n(T'))$.
\end{lemma}
{\it Proof:} 
By definition, $k_1$, $k_2$, and $k_3$ satisfy the restrictions stated in Lemma \ref{lemma3}(\ref{lemma8a}) to (\ref{lemma8h}).
If $k_3\geq 2$, then $(k_1,k_2,k_3)=(0,0,2)$.
Hence, we may assume that $k_3\leq 1$.
We consider different cases.
In what follows, $T'$ will be a tree with $n(T')\leq n(T)-6$,
and $S'$ will be a minimum exponential dominating set of $T'$.
Let $v$ be the parent of $w$.

\medskip

\noindent {\bf Case 1} {\it $k_3=1$.}

\medskip

\noindent In this case $k_2\leq 2$.

If $k_1=0$, then let $T'=T-V(T_w)$, 
and let $Y$ be the set of $2k_2+3$ descendants of $w$ at distance $2$ from $w$.
The set $S'\cup Y$ is an exponential dominating set of $T$.
Since $n(T')=n(T)-5k_2-8$ and $\gamma_e(T)\leq \gamma_e(T')+2k_2+3$,
we obtain $\gamma_e(T)\leq \gamma_e(T')+\frac{7}{18}(n(T)-n(T'))$ as $2k_2+3\le \frac{7}{18}(5k_2+8)$.
Hence, we may assume that $k_1\geq 1$. Now, $k_1\geq 1$ and $k_3=1$ imply $k_1\leq 2$ and $k_2\leq 1$.

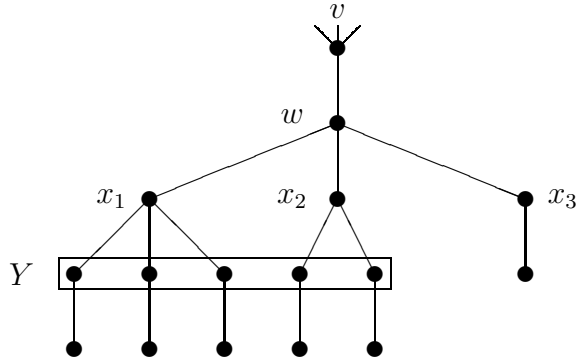
\begin{figure}[H]
\begin{center}
\unitlength 1mm 
\linethickness{0.4pt}
\ifx\plotpoint\undefined\newsavebox{\plotpoint}\fi 
\begin{picture}(75,45)(0,0)
\put(10,0){\circle*{2}}
\put(20,0){\circle*{2}}
\put(30,0){\circle*{2}}
\put(40,0){\circle*{2}}
\put(50,0){\circle*{2}}
\put(10,10){\circle*{2}}
\put(20,10){\circle*{2}}
\put(30,10){\circle*{2}}
\put(40,10){\circle*{2}}
\put(50,10){\circle*{2}}
\put(70,10){\circle*{2}}
\put(10,10){\line(0,-1){10}}
\put(20,10){\line(0,-1){10}}
\put(30,10){\line(0,-1){10}}
\put(40,10){\line(0,-1){10}}
\put(50,10){\line(0,-1){10}}
\put(70,20){\circle*{2}}
\put(70,20){\line(0,-1){10}}
\put(45,40){\circle*{2}}
\put(45,30){\line(5,-2){25}}
\multiput(45,40)(.03370787,.03370787){89}{\line(0,1){.03370787}}
\put(45,40){\line(0,1){3}}
\multiput(45,40)(-.03370787,.03370787){89}{\line(0,1){.03370787}}
\put(42,43){\line(0,1){0}}
\put(45,30){\circle*{2}}
\put(45,20){\circle*{2}}
\put(20,20){\circle*{2}}
\put(50,10){\line(-1,2){5}}
\put(45,20){\line(0,1){20}}
\put(40,10){\line(1,2){5}}
\put(30,10){\line(-1,1){10}}
\put(20,20){\line(5,2){25}}
\put(20,10){\line(0,1){10}}
\put(20,20){\line(-1,-1){10}}
\put(45,45){\makebox(0,0)[cc]{$v$}}
\put(15,20){\makebox(0,0)[cc]{$x_1$}}
\put(39,20){\makebox(0,0)[cc]{$x_2$}}
\put(75,20){\makebox(0,0)[cc]{$x_3$}}
\put(8,8){\framebox(44,4)[cc]{}}
\put(3,10){\makebox(0,0)[cc]{$Y$}}
\put(39,31){\makebox(0,0)[cc]{$w$}}
\end{picture}
\end{center}
\caption{$(k_1,k_2,k_3)=(1,1,1)$.}\label{figlemma9}
\end{figure}

If $k_2=0$, then $(k_1,k_2,k_3)\in \{ (1,0,1),(2,0,1)\}$.
Hence, we may assume that $k_2=1$, which implies $k_1=1$.
Let $x_1$ and $x_2$ be the two children of $w$ of degree at least $3$,
and let $x_3$ be the child of $w$ of degree $2$.
Let $T'=T-(\{ w\}\cup V(T_{x_1})\cup V(T_{x_2}))+vx_3$.
See Figure \ref{figlemma9} for an illustration.

\noindent Note that if $S'$ contains neither $x_3$ nor the child of $x_3$, then $w_{(T',S')}(v)\geq 4$.
Let $Y$ be the set of the five children of $x_1$ and $x_2$.
Since $S'\cup Y$ is an exponential dominating set of $T$,
we obtain $\gamma_e(T)\leq \gamma_e(T')+\frac{5}{13}(n(T)-n(T'))$.

\medskip

\noindent {\bf Case 2} {\it $k_3=0$.}

\medskip

\noindent In this case $k_2\leq 5$.

If $k_1=0$, then let $T'=T-V(T_w)$, and let $Y$ be the set of the $2k_2$ descendants of $w$ at distance $2$ from $w$.
The set $S'\cup Y$ is an exponential dominating set of $T$.
Since $n(T')=n(T)-5k_2-1$ and $\gamma_e(T)\leq \gamma_e(T')+2k_2$,
we obtain $\gamma_e(T)\leq \gamma_e(T')+\frac{5}{13}(n(T)-n(T'))$ as $2k_2\le \frac{5}{13}(5k_2+1)$ for $k_2\le 5$.
Note that $n(T)-n(T')=6$ only for $k_2=1$.
Hence, we may assume that $k_1\geq 1$, which implies $k_2\leq 2$.

\medskip

\noindent {\bf Case 2.1} {\it $k_2=2$.}

\medskip

\noindent In this case $k_1=1$.
Let $x_1$ and $x_2$ be the two children of $w$ of degree $3$,
and let $x_3$ be the child of $w$ of degree $2$.
Let $T'=T-(\{ w\}\cup V(T_{x_1})\cup V(T_{x_2}))+vx_3$.
Let $Y$ be the set of the four children of $x_1$ and $x_2$.
Since $S'\cup Y$ is an exponential dominating set of $T$,
we obtain $\gamma_e(T)\leq \gamma_e(T')+\frac{4}{11}(n(T)-n(T'))$.

\medskip

\noindent {\bf Case 2.2} {\it $k_2=1$.}

\medskip

\noindent In this case $k_1\leq 2$.

If $k_1=2$, then $(k_1,k_2,k_3)=(2,1,0)$.
Hence, we may assume that $k_1=1$.
Let $T'=T-V(T_w)$, and let $X$ be the set containing the three parents of leaves in $T_w$.
Since $S'\cup X$ is an exponential dominating set of $T$,
we obtain $\gamma_e(T)\leq \gamma_e(T')+\frac{3}{8}(n(T)-n(T'))$.
$\Box$

\medskip

\noindent A vertex of $T$ is {\it good} if it has one of the types in $\{ (0,0,2),(1,0,1),(2,0,1),(2,1,0)\}$.

\begin{lemma}\label{lemma5}
If a vertex $v$ of $T$ has two children $w_1$ and $w_2$ such that $w_1$ has type $(2,0,1)$ and $w_2$ is good,
then there is a tree $T'$ with $n(T')<n(T)$ and $\gamma_e(T)\leq \gamma_e(T')+\frac{9}{23}(n(T)-n(T'))$.
\end{lemma}
{\it Proof:} Let $T'$ arise from $T-(V(T_{w_1})\cup V(T_{w_2}))$ by
adding the new vertex $w$, and adding the new edge $vw$.  Let $S'$ be
a minimum exponential dominating set of $T'$.  
For each possible type of $w_2$, we construct an exponential dominating set $S$ of $T$ from
$S'$ as follow: If $v\in S'$, then $S$ is the union of $S'\setminus \{ v\}$
and all parents of leaves of $T_v$, and if $v\not\in S'$, then we
add to $S'$ all parents of leaves of $T_v$ except for one
child of $w_1$ to obtain $S$. We let the reader check that $S$ is an
exponential dominating set of $T$ (using the fact that when $v\not\in
S'$, we must have $w_{(T',S')}(v)\geq 2$) and that we obtain the
following results.\\ If $w_2$ has type $(2,0,1)$, then $n(T')=n(T)-23$
and $\gamma_e(T)\leq \gamma_e(T')+9$.\\ If $w_2$ has type $(1,0,1)$ or
$(2,1,0)$, then $n(T')=n(T)-21$ and $\gamma_e(T)\leq
\gamma_e(T')+8$.\\ If $w_2$ has type $(0,0,2)$, then $n(T')=n(T)-26$
and $\gamma_e(T)\leq \gamma_e(T')+10$.  $\Box$

\begin{lemma}\label{lemma6}
If a vertex of $T$ has three children that are good, 
then there is a tree $T'$ with $n(T')<n(T)$ and $\gamma_e(T)\leq \gamma_e(T')+\frac{13}{33}(n(T)-n(T'))$.
\end{lemma}
{\it Proof:} Suppose that $v$ is a vertex of $T$ that has three good children $w_1$, $w_2$, and $w_3$.
By Lemma \ref{lemma5}, no child of $v$ has type $(2,0,1)$.
For each $(k_1,k_2,k_3)\in \{ (1,0,1),(2,1,0),(0,0,2)\}$, let $n(k_1,k_2,k_3)$ vertices in $\{ w_1,w_2,w_3\}$ have type $(k_1,k_2,k_3)$.

First, suppose that $n(1,0,1)\geq 2$ and say that $w_1$ and $w_2$ are
of type $(1,0,1)$.  Let $T'$ arise from $T-(V(T_{w_1})\cup
V(T_{w_2})\cup V(T_{w_3}))$ by adding the two new vertices $w$ and
$x$, and adding the two new edges $vw$ and $wx$.  
From an exponential dominating set $S'$ of $T'$, we construct an exponential dominating set of
$T$.  If $w$ or $x$ belongs to $S'$, then $S$ is the union of
$S'\setminus \{x,w\}$ and all the leaves of $T_{w_1}$ and $T_{w_2}$. If
$w$ and $x$ do not belong to $S'$, then $v$ does not belong to $S'$
also and $S$ is obtained as the union of $S'$ and all the leaves of
$T_{w_1}$ and $T_{w_2}$ except for one child of $w_1$.  Notice that in
this latter case we have $w_{(T',S')}(v)\ge 4$. In both cases we
obtain the following results.\\ If $n(0,0,2)=0$, then
$n(T')=n(T)-28$ and $\gamma_e(T)\leq \gamma_e(T')+11$.\\ If
$n(0,0,2)=1$, then $n(T')=n(T)-33$ and $\gamma_e(T)\leq
\gamma_e(T')+13$.\\ Hence, in these cases $\gamma_e(T)\leq
\gamma_e(T')+\frac{13}{33}(n(T)-n(T'))$.

Next, suppose that either $n(1,0,1)=1$ or $n(0,0,2)=3$.
Let $T'$ arise from $T-(V(T_{w_1})\cup V(T_{w_2})\cup V(T_{w_3}))$
by adding the new vertex $w$, and adding the new edge $vw$.
We derive as previously an exponential dominating set of $T$ from
an exponential dominating set of $T'$ and obtain the following results.\\
If $n(1,0,1)=1$, then considering the three possibilities for the other values, 
it follows that $\gamma_e(T)\leq \gamma_e(T')+\frac{15}{39}(n(T)-n(T'))$.\\
If $n(0,0,2)=3$,
then $n(T')=n(T)-44$ and $\gamma_e(T)\leq \gamma_e(T')+17$.\\
Hence, in these cases $\gamma_e(T)\leq \gamma_e(T')+\frac{17}{44}(n(T)-n(T'))$.

In what follows, we may assume that $n(1,0,1)=0$ and $n(0,0,2)\leq 2$.
Let $T'=T-(V(T_{w_1})\cup V(T_{w_2})\cup V(T_{w_3}))$. Once again we define an 
exponential dominating set of $T$ from one of $T'$ and obtain the following.
We have $n(T')=n(T)-10n(2,1,0)-15n(0,0,2)$ and 
since $n(2,1,0)\geq 1$, it follows that $\gamma_e(T)\leq \gamma_e(T')+4n(2,1,0)+6n(0,0,2)-1$.
Considering the three possibilities for the value of $n(0,0,2)$ 
implies $\gamma_e(T)\leq \gamma_e(T')+\frac{15}{40}(n(T)-n(T'))$.
$\Box$

\medskip

\noindent For the rest of this subsection, let $v$ be a vertex of $T$ such that 
\begin{itemize}
\item $T_v$ has depth $4$,
\item $v$ has at most two children $w$ such that $T_w$ has depth $3$,
\item every child $w$ of $v$ such that $T_w$ has depth $3$ is good, and
\item if $v$ has two children that are good, then none of the two has type $(2,0,1)$.
\end{itemize}
Let $W$ be the set of children $w$ of $v$ such that $T_w$ has depth $3$.
Let $T^{(0)}=T-\bigcup_{w\in W}V(T_w)$, and let $d_{\rm red}$ be the depth of $T^{(0)}_v$.
By construction, $d_{\rm red}\leq 3$.
Note that, since $T$ has depth at least $5$, the vertex $v$ has a parent $u$ in $T$.

\begin{lemma}\label{lemma8}
If $d_{\rm red}\leq 2$, 
then there is a tree $T'$ with $n(T')<n(T)$ and $\gamma_e(T)\leq \gamma_e(T')+\frac{13}{33}(n(T)-n(T'))$.
\end{lemma}
{\it Proof:} First, suppose that $d_{\rm red}=0$.
Let $T'=T-V(T_v)$.
If $|W|=1$, then 
$\gamma_e(T)\leq \gamma_e(T')+\frac{5}{13}(n(T)-n(T'))$.
If $|W|=2$, then  
$\gamma_e(T)\leq \gamma_e(T')+\frac{12}{31}(n(T)-n(T'))$.
Next, suppose that $d_{\rm red}=1$.
Let $T'$ arise from $T$ by removing all descendants of $v$.
For the two following cases we simply extend an exponential dominating set of
$T'$ by adding all the parents of the leaves of $T_v$.
If $|W|=1$, then 
$\gamma_e(T)\leq \gamma_e(T')+\frac{5}{13}(n(T)-n(T'))$.
If $|W|=2$, then 
$\gamma_e(T)\leq \gamma_e(T')+\frac{12}{31}(n(T)-n(T'))$.

Hence, we may assume that $d_{\rm red}=2$.  Let $v$ have $n_1$
children that are leaves, and $n_2$ children $w$ such that $T_w$ has
depth $1$.  Since $d_{\rm red}=2$, we have $n_2\geq 1$.  By Lemma
\ref{lemma2}(\ref{r1}), we may assume that $T_w\cong T(1,0)$ for every
child $w$ of $v$ such that $T_w$ has depth $1$.  We argue as
previously for the following cases.

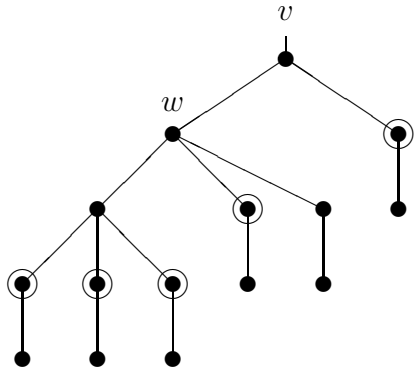
\begin{figure}[H]
\begin{center}
\unitlength 1mm 
\linethickness{0.4pt}
\ifx\plotpoint\undefined\newsavebox{\plotpoint}\fi 
\begin{picture}(52,46)(0,0)
\put(0,0){\circle*{2}}
\put(10,0){\circle*{2}}
\put(20,0){\circle*{2}}
\put(0,10){\circle*{2}}
\put(10,10){\circle*{2}}
\put(20,10){\circle*{2}}
\put(0,10){\line(0,-1){10}}
\put(10,10){\line(0,-1){10}}
\put(20,10){\line(0,-1){10}}
\put(10,20){\circle*{2}}
\put(20,10){\line(-1,1){10}}
\put(10,10){\line(0,1){10}}
\put(10,20){\line(-1,-1){10}}
\put(30,10){\circle*{2}}
\put(40,10){\circle*{2}}
\put(30,20){\circle*{2}}
\put(40,20){\circle*{2}}
\put(40,10){\line(0,1){10}}
\put(30,20){\line(0,-1){10}}
\put(20,30){\circle*{2}}
\put(50,20){\circle*{2}}
\put(50,30){\circle*{2}}
\put(35,40){\circle*{2}}
\put(50,20){\line(0,1){10}}
\put(50,30){\line(-3,2){15}}
\put(35,40){\line(-3,-2){15}}
\put(34,40){\line(1,0){1}}
\put(35,40){\line(0,1){0}}
\put(35,43){\line(0,-1){3}}
\put(35,46){\makebox(0,0)[cc]{$v$}}
\put(20,34){\makebox(0,0)[cc]{$w$}}
\put(0,10){\circle{4}}
\put(10,10){\circle{4}}
\put(20,10){\circle{4}}
\put(30,20){\circle{4}}
\put(50,30){\circle{4}}
\put(10,20){\line(1,1){10}}
\put(20,30){\line(1,-1){10}}
\put(40,20){\line(-2,1){20}}
\end{picture}
\end{center}
\caption{$W=\{ w\}$, $w$ has type $(2,0,1)$, $n_1=0$, and $n_2=1$.}\label{figlemma12}
\end{figure}

First, suppose that $n_2=1$.  If $v$ has a child $w$ of type
$(2,0,1)$, then $w$ is the unique child of $v$ such that $T_w$ has
depth $3$.  See Figure \ref{figlemma12} for an illustration.

\noindent In this case, let $T'$ arise from $T$ by removing all descendants of $v$ except for $w$.
Now, $n(T')=n(T)-13-n_1$ and $\gamma_e(T)\leq \gamma_e(T')+5$,
which implies $\gamma_e(T)\leq \gamma_e(T')+\frac{5}{13}(n(T)-n(T'))$.
Hence, we may assume that $v$ has no child of type $(2,0,1)$.
Now, for $T'=T-V(T_v)$, it follows that $\gamma_e(T)\leq \gamma_e(T')+\frac{13}{33}(n(T)-n(T'))$.

Next, suppose that $n_2\geq 3$.
Let $T'$ arise from $T$ by removing all descendants of $v$.
If $v$ has exactly one child $w$ such that $T_w$ has depth $3$, 
then $\gamma_e(T)\leq \gamma_e(T')+\frac{7}{18}(n(T)-n(T'))$.
If $v$ has two children $w$ such that $T_w$ has depth $3$, 
then $\gamma_e(T)\leq \gamma_e(T')+\frac{14}{36}(n(T)-n(T'))$.

Finally, suppose that $n_2=2$.
If $v$ has a child $w$ of type $(2,0,1)$, then $w$ is the unique child of $v$ such that $T_w$ has depth $3$.
In this case, let $T'$ arise from $T$ by removing all descendants of $v$.
Now, $n(T')=n(T)-16-n_1$ and $\gamma_e(T)\leq \gamma_e(T')+6$,
which implies $\gamma_e(T)\leq \gamma_e(T')+\frac{3}{8}(n(T)-n(T'))$.
Hence, we may assume that $v$ has no child of type $(2,0,1)$.
Let $T'$ arise from $T$ by removing all descendants of $v$ except for one child of $v$.
It follows that $\gamma_e(T)\leq \gamma_e(T')+\frac{13}{33}(n(T)-n(T'))$.
$\Box$

\subsection{Proof of Theorem \ref{theorem1}}\label{sec2.2}

Since $\gamma_e(G)\leq \gamma_e(H)$ for every spanning subgraph $H$ of
$G$, it suffices to prove the statement in the case that $G$ is a tree
$T$.  For a contradiction, suppose that $T$ is a counterexample of
minimum order.  This choice of $T$ implies that there is no tree $T'$
with $n(T')<n(T)$ and $\gamma_e(T)\leq
\gamma_e(T')+\alpha(n(T)-n(T'))$ for some $\alpha\leq \frac{43}{108}$.
By Lemma \ref{lemma1}, $T$ has diameter at least $5$.  Root $T$ in a
vertex of maximum eccentricity.  Let $v$ be a vertex of $T$ such that
$T_v$ has depth $4$.  Let $W$ be the set of children $w$ of $v$ such
that $T_w$ has depth $3$.  By Lemma \ref{lemma3} and Lemma
\ref{lemma4}, every vertex in $W$ is good.  By Lemma \ref{lemma6},
$|W|\leq 2$, and, by Lemma \ref{lemma5}, if $|W|=2$, then no vertex in
$W$ has type $(2,0,1)$.  Let $T^{(0)}=T-\bigcup_{w\in W}V(T_w)$, and
let $d_{\rm red}$ be the depth of $T^{(0)}_v$.  By construction,
$d_{\rm red}\leq 3$, and Lemma \ref{lemma8} implies $d_{\rm red}=3$.

Now let $X$ be the set of children $x$ of $v$ such that $T_x$ has
depth at most $2$.  By Lemma \ref{lemma3} applied to $v$ and $X$, the
vertex $v$ has some type in the rooted tree $T^{(0)}$.  Let
$T^{(1)}=T-V(T_v)$.  Let $T^{(2)}$ arise from $T$ by removing all
descendants of $v$.  Finally, let $T^{(3)}$ arise from $T$ by removing
all descendants of $v$ except for one child of $v$.  As before we will
extend a minimum exponential dominating set $S^{(i)}$ of some
$T^{(i)}$ to obtain an exponential dominating set of $T$. 
We will use that $w_{(T^{(1)},S^{(1)})}(u)\ge 1$
where $u$ is the parent of $v$ in $T$, that
$w_{(T^{(2)},S^{(2)})}(v)\ge 1$, 
and that $w_{(T^{(3)},S^{(3)})}(v)\ge 2$ assuming that the child of $v$ in $T^{(3)}$ does not belong to
$S^{(3)}$. 
As before, for computations with $T^{(3)}$, 
we have to distinguish the cases $v\in S^{(3)}$ and $v\notin S^{(3)}$.

First, suppose that no vertex in $W$ has type $(2,0,1)$.
This implies that $n(T)-n(T^{(0)})\leq 30$ and 
$\gamma_e(T)\leq \gamma_e(T^{(0)})+\frac{2}{5}(n(T)-n(T^{(0)}))$
(by taking all the parents of the leaves of $T_w$ for $w\in W$ to extend
an exponential dominating set of $T^{(0)}$).
If $v$ is not good in $T^{(0)}$, 
then Lemma \ref{lemma4} implies that there is a tree $T'$
with $n(T^{(0)})-n(T')\geq 6$ and $\gamma_e(T^{(0)})\leq \gamma_e(T')+\frac{7}{18}(n(T^{(0)})-n(T'))$.
If $0<\alpha_1<\alpha_2$, $1\leq n_1^0\leq n_1$, and $0\leq n_2\leq n_2^0$, then
$\alpha_1n_1+\alpha_2n_2\leq (\frac{\alpha_1n^0_1+\alpha_2n^0_2}{n^0_1+n^0_2})(n_1+n_2)$.
Therefore,
\begin{eqnarray*}
\gamma_e(T) & \leq & \gamma_e(T^{(0)})+\frac{2}{5}(n(T)-n(T^{(0)}))\\
& \leq & \gamma_e(T')+\frac{7}{18}(n(T^{(0)})-n(T'))+\frac{2}{5}(n(T)-n(T^{(0)}))\\
& \leq & \gamma_e(T')+\frac{\frac{7}{18}\cdot 6+\frac{2}{5}\cdot 30}{6+30}(n(T)-n(T'))\\
& \leq & \gamma_e(T')+\frac{43}{108}(n(T)-n(T')),
\end{eqnarray*}
which is a contradiction.
Hence, $v$ is good in $T^{(0)}$.
If $v$ is of type $(2,0,1)$ in $T^{(0)}$, 
then adding all parents of leaves in $T_v$ except for one child of $v$ 
to a minimum exponential dominating set of $T^{(2)}$ that does not contain $v$
yields an exponential dominating set of $T$.
This implies 
$\gamma_e(T)\leq \gamma_e(T^{(2)})+\frac{16}{41}(n(T)-n(T^{(2)}))$,
which is a contradiction (the worst case appears when $|W|=2$ and each $w\in W$ is of type $(0,0,2)$).
Hence, we may assume that $v$ is not of type $(2,0,1)$ in $T^{(0)}$.
It follows that 
$\gamma_e(T)\leq \gamma_e(T^{(3)})+\frac{17}{43}(n(T)-n(T^{(3)}))$,
which is a contradiction (the worst case appears when $|W|=2$ and each $w\in W$ is of type $(0,0,2)$ and $v$ is of 
type $(0,0,2)$ in $T^{(0)}$).

Hence, we may assume that one child $w$ of $v$ is of type $(2,0,1)$, which implies that $w$ is the only element of $W$ by
Lemma~\ref{lemma5}. Let $v$ have type $(k_1,k_2,k_3)$ in $T^{(0)}$.

First, suppose $k_3=2$.  This implies $(k_1,k_2,k_3)=(0,0,2)$ by
Lemma~\ref{lemma3}(\ref{lemma8a}), $n(T^{(1)})=n(T)-27$, and
$\gamma_e(T)\leq \gamma_e(T^{(1)})+10$ by adding to a minimum
exponential dominating set of $T^{(1)}$ all the parents of the leaves
of $T_v$ except for one child of $w$. So we have a contradiction and
we assume that $k_3\leq 1$.

\begin{figure}[H]
\begin{center}
\unitlength 1mm 
\linethickness{0.4pt}
\ifx\plotpoint\undefined\newsavebox{\plotpoint}\fi 
\begin{picture}(101,44)(0,0)
\put(0,0){\circle*{2}}
\put(50,10){\circle*{2}}
\put(10,0){\circle*{2}}
\put(60,10){\circle*{2}}
\put(80,10){\circle*{2}}
\put(20,0){\circle*{2}}
\put(70,10){\circle*{2}}
\put(90,10){\circle*{2}}
\put(0,10){\circle*{2}}
\put(50,20){\circle*{2}}
\put(10,10){\circle*{2}}
\put(60,20){\circle*{2}}
\put(80,20){\circle*{2}}
\put(20,10){\circle*{2}}
\put(70,20){\circle*{2}}
\put(90,20){\circle*{2}}
\put(0,10){\line(0,-1){10}}
\put(50,20){\line(0,-1){10}}
\put(10,10){\line(0,-1){10}}
\put(60,20){\line(0,-1){10}}
\put(80,20){\line(0,-1){10}}
\put(20,10){\line(0,-1){10}}
\put(70,20){\line(0,-1){10}}
\put(90,20){\line(0,-1){10}}
\put(10,20){\circle*{2}}
\put(20,10){\line(-1,1){10}}
\put(10,10){\line(0,1){10}}
\put(10,20){\line(-1,-1){10}}
\put(30,10){\circle*{2}}
\put(40,10){\circle*{2}}
\put(30,20){\circle*{2}}
\put(40,20){\circle*{2}}
\put(40,10){\line(0,1){10}}
\put(30,20){\line(0,-1){10}}
\put(20,30){\circle*{2}}
\put(20,34){\makebox(0,0)[cc]{$w$}}
\put(0,10){\circle{4}}
\put(50,20){\circle{4}}
\put(10,10){\circle{4}}
\put(60,20){\circle{4}}
\put(80,20){\circle{4}}
\put(20,10){\circle{4}}
\put(70,20){\circle{4}}
\put(90,20){\circle{4}}
\put(30,20){\circle{4}}
\put(10,20){\line(1,1){10}}
\put(20,30){\line(1,-1){10}}
\put(40,20){\line(-2,1){20}}
\put(60,30){\circle*{2}}
\put(85,30){\circle*{2}}
\put(90,20){\line(-1,2){5}}
\put(85,30){\line(-1,-2){5}}
\put(70,20){\line(-1,1){10}}
\put(60,30){\line(0,-1){10}}
\put(50,20){\line(1,1){10}}
\put(50,40){\circle*{2}}
\put(20,30){\line(3,1){30}}
\put(50,40){\line(1,-1){10}}
\multiput(50,40)(.1178451178,-.0336700337){297}{\line(1,0){.1178451178}}
\put(100,20){\circle*{2}}
\put(100,30){\circle*{2}}
\put(50,40){\line(5,-1){50}}
\put(100,30){\line(0,-1){10}}
\put(50,40){\line(0,1){4}}
\put(46,42){\makebox(0,0)[cc]{$v$}}
\end{picture}
\end{center}
\caption{$w$ has type $(2,0,1)$ in $T$, and $v$ has type $(k_1,k_2,k_3)=(1,1,1)$ in $T^{(0)}$.}\label{fig12}
\end{figure}
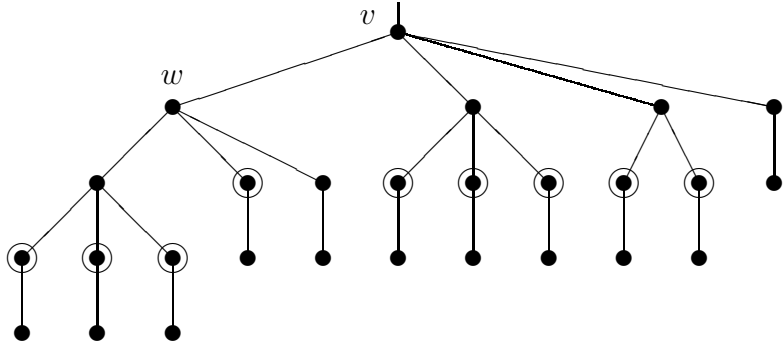

Next, suppose $k_3=1$.
In this case Lemma~\ref{lemma3}(\ref{lemma8b}) implies $k_2\leq 2$.
If $k_2=1$, 
then adding all parents of leaves in $T_v$ 
that are no children of $v$ 
except for one child of $w$ 
to a minimum exponential dominating set of $T^{(1)}$
yields an exponential dominating set of $T$.
See Figure \ref{fig12} for an illustration.

\noindent This implies $n(T^{(1)})\leq n(T)-25-2k_1$ and $\gamma_e(T)\leq \gamma_e(T^{(1)})+9$,
which is a contradiction.
Similarly, if $k_2=2$, then 
$n(T^{(1)})\leq n(T)-30-2k_1$ and $\gamma_e(T)\leq \gamma_e(T^{(1)})+11$,
which is a contradiction.
Hence, we may assume that $k_2=0$.
By Lemma \ref{lemma2}(\ref{r14}), we have $k_1\leq 2$.
If $k_1=1$, then
$n(T^{(1)})\leq n(T)-22$ and $\gamma_e(T)\leq \gamma_e(T^{(1)})+8$,
which is a contradiction.
If $k_1=2$, then
$n(T^{(1)})\leq n(T)-24$ and $\gamma_e(T)\leq \gamma_e(T^{(1)})+9$,
which is a contradiction.
Hence, we may assume that $k_1=0$.
Now $n(T^{(3)})=n(T)-18$ and $\gamma_e(T)\leq \gamma_e(T^{(3)})+7$,
which is also a contradiction.
Hence, we may assume that $k_3=0$.

In this case Lemma~\ref{lemma3}(\ref{lemma8e}) implies $1\leq k_2\leq 5$.
If $k_1=0$, then we have $n(T^{(3)})\leq n(T)-5k_2-11$ and $\gamma_e(T)\leq
\gamma_e(T^{(3)})+2k_2+4$ implying the contradiction $\gamma_e(T)\leq
\gamma_e(T^{(3)})+\frac{14}{36}(n(T)-n(T^{(3)}))$. 
Hence, we may assume that $k_1\geq 1$.
By Lemma~\ref{lemma3}(\ref{lemma8f}), we have $k_2\leq 2$.  By
Lemma~\ref{lemma3}(\ref{lemma8g}), if $k_2=2$, then we have $k_1=1$.  
In this case we obtain $n(T^{(3)})\leq n(T)-23$ and $\gamma_e(T)\leq
\gamma_e(T^{(3)})+9$, a contradiction.  
Finally, if $k_2=1$, then, by Lemma~\ref{lemma3}(\ref{lemma8h}), we have $k_1\le 2$.  
If $k_1=1$, then $n(T^{(3)})\leq n(T)-18$ and $\gamma_e(T)\leq
\gamma_e(T^{(3)})+7$, a contradiction. 
If $k_1=2$, then $n(T^{(2)})\leq n(T)-21$ and $\gamma_e(T)\leq \gamma_e(T^{(2)})+8$
implying the contradiction $\gamma_e(T)\leq
\gamma_e(T^{(2)})+\frac{8}{21}(n(T)-n(T^{(2)}))$, which completes the
proof.  $\Box$


\begin{thebibliography}{}
\bibitem{adr} J.D. Alvarado, S. Dantas, D. Rautenbach, Distance $k$-domination, distance $k$-guarding, and distance $k$-vertex cover of maximal outerplanar graphs, Discrete Appl. Math. 194 (2015) 154-159.
\bibitem{bz} G. Bacs\'{o}, Zs. Tuza, Distance domination versus iterated domination, Discrete Math. 312 (2012) 2672-2675.
\bibitem{chmm} S. Canales, G. Hern\'{a}ndez, M. Martins, I. Matos, Distance domination, guarding and covering of maximal outerplanar graphs, Discrete Appl. Math. 181 (2015) 41-49.
\bibitem{cr} Y. Caro, Y. Roditty, A note on the $k$-domination number of a graph, Int. J. Math. Math. Sci. 13 (1990) 205-206.
\bibitem{cgs} E.J. Cockayne, B. Gamble, B. Shepherd, An upper bound for the $k$-domination number of a graph, J. Graph Theory 9 (1985) 533-534.
\bibitem{chm} E.J. Cockayne, S. Herke, C.M. Mynhardt, Broadcasts and domination in trees, Discrete Math. 311 (2011) 1235-1246.
\bibitem{ddems} P. Dankelmann, D. Day, D. Erwin, S. Mukwembi, H. Swart, Domination with exponential decay, Discrete Math. 309 (2009) 5877-5883.
\bibitem{dghpv} E. DeLaVi\~{n}a, W. Goddard, M.A. Henning, R. Pepper, E.R. Vaughan, Bounds on the $k$-domination number of a graph, Appl. Math. Lett. 24 (2011) 996-998.
\bibitem{dehhh} J.E. Dunbar, D.J. Erwin, T.W. Haynes, S.M. Hedetniemi, S.T. Hedetniemi, Broadcasts in graphs, Discrete App. Math. 154 (2006) 59-75.
\bibitem{e} D.J. Erwin, Dominating broadcasts in graphs, Bull. Inst. Comb. Appl. 42 (2004) 89-105.
\bibitem{fhv} O. Favaron, A. Hansberg, L. Volkmann, On $k$-domination and minimum degree in graphs, J. Graph Theory 57 (2008) 33-40.
\bibitem{fj} J.F. Fink, M.S. Jacobson, $n$-domination in graphs, in: Graph Theory with Applications to Algorithms and Computer Science, John Wiley and Sons, 1985, 282-300.
\bibitem{hmv} A. Hansberg, D. Meierling, L. Volkmann, Distance domination and distance irredundance in graphs, Electron. J. Combin. 14 (2007), $\#$R35.
\bibitem{ha} A. Hansberg, Bounds on the connected $k$-domination number in graphs, Discrete Appl. Math. 158 (2010) 1506-1510.
\bibitem{hv} A. Hansberg, L. Volkmann, Upper bounds on the $k$-domination number and the $k$-Roman domination number, Discrete Appl. Math. 157 (2009) 1634-1639.
\bibitem{hhs} T.W. Haynes, S.T. Hedetniemi, P.J. Slater, Fundamentals of Domination in Graphs, Marcel Dekker, Inc., New York, 1998.
\bibitem{he} M.A. Henning, Distance domination in graphs, in: T.W. Haynes, S.T. Hedetniemi, P.J. Slater (Eds.), Domination in Graphs: Advanced Topics, Marcel Dekker, New York, 1998, 321-349.
\bibitem{hl} P. Heggernes, D. Lokshtanov, Optimal broadcast domination in polynomial time, Discrete Math. 306 (2006) 3267-3280.
\bibitem{rv} D. Rautenbach, L. Volkmann, New bounds on the $k$-domination number and the $k$-tuple domination number, Appl. Math. Lett. 20 (2007) 98-102.
\bibitem{tx} F. Tian, J.-M. Xu, A note on distance domination numbers of graphs, Australas. J. Combin. 43 (2009) 181-190.
\end{thebibliography}
\end{document}